\documentclass{amsart}

\usepackage{amsmath, amssymb,epsfig, verbatim}
\usepackage{graphicx}
\usepackage{epic}
\usepackage{color}

\newtheorem{theorem}{Theorem}[section]

\newtheorem{proposition}[theorem]{Proposition}

\newtheorem{remark}[theorem]{Remark}
\numberwithin{equation}{section}

\newcommand{\ZZ}{\mathbf Z}

\newcommand{\FF}{\mathbf F}
\newcommand{\HH}{\mathbf H}

\newcommand{\CP}{\mathbf{CP}}
\newcommand{\RP}{\mathbf{RP}}

\begin{document}

\title{A fake quadric that is not a bidisk quotient}
\author{Zolt\'an Szab\'o}
\address{Department of Mathematics, Princeton University,
Princeton, NJ 08544}
\email{szabo@math.princeton.edu}
\date{November 24, 2026}

\begin{abstract}
  Fake quadrics are minimal complex surfaces of general type with $b_1=0$ and $e=4$. The paper gives a construction of a
  fake quadric $S$ that is a $32$-fold branched cover of $\mathbf{CP}^1\times\mathbf{CP}^1$ along nine curves. This configuration is given by four horizontal lines, four vertical lines and the diagonal.
  We show that   $S$ is not aspherical. As a corollary
  the universal cover of $S$ is not the bidisk $\mathbf{H}\times\mathbf{H}$.

\end{abstract}

\maketitle

\section{Introduction}

A fake quadric is a smooth minimal compact complex surface of general type that has  the same rational homology
groups as $\CP^1\times\CP^1$. Many examples are quotients of the bidisk
$\HH\times\HH$ and hence are aspherical. These include surfaces that admit a finite unramified cover
which is  a product of two  curves \cite{Cat2},  and irreducible  arithmetic examples
\cite{Dzambic,LSV,Shavel}. In a different direction
Rito \cite{Rito} and Borisov-Rito \cite{BR} have constructed fake
quadrics obtained as covers of singular Godeaux surfaces. It is a natural question whether  every fake quadric is
uniformized by the bidisk, see for example \cite[Question~6]{BCP},  \cite[Problem~5]{CCCK}.
In this paper we give an explicit construction of a fake quadric that is not aspherical.
The construction starts with nine curves in $\CP^1\times\CP^1$.

Let $Q=\CP^1\times\CP^1$ and $T=\{1,2,3,4\}$. For each $p\in T$, take the vertical line $\{p\}\times\CP^1$ and the horizontal line
$\CP^1\times\{p\}$. These eight lines and the
diagonal
form a divisor $D$ with four triple points and twelve double points. We consider
branched covers of $Q$ of degree $2^n$ with branch locus $D$.

Fix the basepoint $P_0=(i,-i)$. For $p\in T$, let $x_p$ be the loop in $\CP^1\times\{-i\}$ that follows
a straight arc toward $(p,-i)$, goes once counterclockwise around a small circle about that point  and returns along the same arc.
Define $y_p$ in the same way in $\{i\}\times\CP^1$, using a small circle about $(i,p)$. With these generators, the fundamental group
$ \pi_1(Q\setminus D,P_0)$ has the presentation:
$$ \langle  x_1,x_2,x_3,x_4,\  y_1,y_2,y_3,y_4 \ | \ [x_p,y_q]=1 \ \ {\rm for }\ p\ne q,\ \ x_1x_2x_3x_4=y_4y_3y_2y_1 \rangle .$$

It will be convenient to label the $32$ sheets by $v=(a,b,c,d,e)\in\FF_2^5$. For each generator $\gamma$ write its monodromy as
$$ f_\gamma(v)=v+G_\gamma(v). $$
It suffices to give the eight maps $G_{x_p}$ and $G_{y_p}$.
\begin{align*}
G_{x_1}&=(1,0,1,b,d+ab), & G_{x_2}&=(0,0,1,0,d+ab),\\
G_{x_3}&=(1,1,0,0,0),    & G_{x_4}&=(0,1,0,0,0),\\
G_{y_1}&=(0,0,1,1,c+d),  & G_{y_2}&=(1,0,0,0,d),\\
G_{y_3}&=(0,1,0,a,1),    & G_{y_4}&=(1,1,1,a+b,c).
\end{align*}

We traverse products of loops from left to right. With this convention the permutations $f_{x_p}$ and $f_{y_p}$ satisfy the relations of the presentation and act
transitively on $\FF_2^5$. They define a connected $32$-fold branched cover $S\to Q$.

A positive meridian of
the diagonal is represented by
$$ \delta=x_4^{-1}x_3^{-1}x_2^{-1}x_1^{-1} =y_1^{-1}y_2^{-1}y_3^{-1}y_4^{-1}, $$
and its monodromy is given by
$$ G_\delta=(0,0,0,b,0). $$

Forgetting the last one or two digits gives a tower of covers $S\to X\to Z\to Q$. The cover $Z$, given by the first three digits, is a $(\ZZ/2)^3$-cover branched only along the four vertical and four horizontal lines. Note that its minimal resolution is an Enriques surface (Remark~\ref{remark1}).
The $16$-fold cover $X$, given by the first four digits, is studied in Section 2.

Note also that
the  formulas for $G_{x_p}$ and $G_{y_p}$ are independent of $e$. It follows that changing the last digit from $e$ to $e+1$
commutes with the monodromy and defines an involution $\tau$ on the branched cover $S$. We shall see that  $S$ is a smooth fake quadric and that its quotient $X=S/\tau$ has six $A_1$ singularities. Resolving these six points gives a minimal surface $Y$ of general type with $K_Y^2=2$, $p_g(Y)=q(Y)=0$, and $H_1(Y;\ZZ)\cong\ZZ/2\oplus\ZZ/4$.

The main result is the following:

\begin{theorem}
  Let $S\to Q$ be the $32$-fold branched cover determined by the monodromy above. Then $S$ is a
  smooth minimal compact complex surface of general type with Euler characteristic $e(S)=4$ and first homology group
  $$ H_1(S;\ZZ)\cong(\ZZ/2)^3\oplus\ZZ/4. $$
  In particular, $S$ is a fake quadric. Furthermore it is not aspherical,
  so its universal cover is not biholomorphic to $\HH\times\HH$.

\end{theorem}

In Section 2 we look at the $16$-fold cover and in Section 3 the $32$-fold cover.
In Section 4 we show that $S$ is not aspherical. The main idea is to give a map $\RP^2\to S$ that is
nontrivial on $\pi_1$.

{\it Acknowledgments.}
The author would like to thank
Inan\c{c} Baykur, Peter Ozsv\'ath and Andr\'as Stipsicz for helpful discussions. This work was
partially supported by the Simons Collaboration Grant on New
Structures in Low-Dimensional Topology.

{\it Statement on the use of AI.}
There are Python programs attached
to the arXiv version of this paper. These contain the calculations for the
first homology groups of $X$ and $S$. The AI assistants ChatGPT and Claude were used to implement these algorithms and also  do preliminary calculations for other examples.  These LLMs  were also used for typesetting and language
editing.
The author assumes intellectual responsibility for the
mathematical statements and for any remaining errors.

\section{The degree-$16$ cover}

Let $U=\CP^1\setminus T$. Then $Q\setminus D$ is the space of ordered pairs of distinct points of $U$. Projection onto the first
coordinate is a  fibration with fiber
$U\setminus\{i\}$ over $i$. The loops $y_1,\ldots,y_4$ freely generate the fundamental
group of this fiber, and the loops $x_1,\ldots,x_4$ project to the four meridians in $U$. The presentation of $pi_1$ is straightforward.

Now let's look at the monodromies. Forgetting $e$ gives the degree-$16$ cover $X\to Q$. On its sixteen
sheets each $x_p$ and $y_p$ acts by eight transpositions. The diagonal
meridian sends $d$ to $d+b$ and fixes the other digits, so it has
eight fixed points and four transpositions.

At a double point $(p,q)$ one can choose  local coordinates $u,v$ in which
$D=\{uv=0\}$. A four-sheeted orbit gives the smooth cover
$(s,t)\mapsto(s^2,t^2)$. A two-sheeted orbit gives the double cover
$w^2=uv$, with an $A_1$ singularity at the origin.

Two-sheeted orbits occur only over $(4,3)$, where $x_4$ and $y_3$ agree  on the sheets with $a=0$. These
sheets give four $A_1$ points, one for each choice of $(c,d)$. The remaining eight
sheets give two smooth points.

Over $(1,1)$ the three monodromies are $x_1$ around the vertical line,
$y_1$ around the horizontal line, and $\delta$ around the diagonal.
They fix $b$ and act on the other three digits by
$$ x_1(a,c,d)=(a+1,c+1,d+b),\ y_1(a,c,d)=(a,c+1,d+1),\  \delta(a,c,d)=(a,c,d+b). $$
When $b=0$, the diagonal monodromy is trivial and both points are smooth.
When $b=1$, the eight sheets give one $A_1$ singularity. This can be seen as follows.
Forgetting $d$ gives the smooth cover $(s,t)\mapsto(s^2,t^2)$, in which the diagonal becomes the two lines $s=\pm t$. Then  the  last digit $d$
gives the double cover $w^2=s^2-t^2$, branched along both lines. 

At $(p,p)$ the three monodromies are $x_p$, $y_p$, and the local
diagonal meridian
$$ \delta_p=(y_{p-1}\cdots y_1)\,\delta\,(y_{p-1}\cdots y_1)^{-1},\ \  \delta_1=\delta. $$
On the sixteen sheets, $\delta_2,\delta_3,\delta_4$ fix $a,b,c$ and act
on $d$ by
$$ \delta_2\colon d\mapsto d+b,\ \ \delta_3\colon d\mapsto d+b,\ \ \delta_4\colon d\mapsto d+b+1.$$
Over $(2,2)$ there are again two smooth points with $b=0$ and one
$A_1$ point with $b=1$. Over each of $(3,3)$ and $(4,4)$ there are
two smooth points. Thus the triple points contribute two $A_1$
singularities, one over $(1,1)$ and one over $(2,2)$.

An easy Euler characteristic computation gives $e(X)=4$.
Let $C_1$ and $C_2$ be the preimages in $X$ of a general vertical and a general horizontal line. Each is a $16$-fold cover of a line
with $36$ simple ramification points, so $C_1$ and $C_2$ have genus three.

\begin{remark}\label{remark1}
  Lifting the presentation of $\pi_1(Q\setminus D)$ to the sixteen sheets, filling in the lifted meridians and
  killing the fundamental groups of the links of the nodes gives $H_1(X;\ZZ)\cong\ZZ/2\oplus\ZZ/4$.
  Adjunction on $C_1$ and $C_2$ gives that over rational coefficients $K_X$ is equal to $(C_1+C_2)/4$,
  so the minimal resolution $Y$ of $X$ is a minimal surface of general type with $K_Y^2=2$, $p_g=q=0$ and the same $H_1$.
  Forgetting $d$ gives the $(\ZZ/2)^3$-cover $Z\to Q$ branched along the four vertical and four horizontal lines.
  Its minimal resolution has an an elliptic fibration with Euler characteristic $12$  and two fibers with multipicity $2$.
\end{remark}

\section{The degree-$32$ cover}\label{sec:S}

Let $z$ and $w$ be the coordinates on the two factors of $Q$. Forgetting the last digit $e$
gives a double cover $f\colon S\to X$. The involution $\tau$ changes $e$ to $e+1$. Let $B$ be the curve in
$X$ over the line $w=4$. It is smooth of genus two and misses the six nodes. Applying $y_4$ twice fixes the first four digits and changes
$e$, so $y_4^2=\tau$.

On the $32$ sheets, $x_1,\dots,x_4$, $y_1,y_2,y_3$ and $\delta$ are still involutions, so $f$ is unramified
at the smooth points of $X$ off $B$. Along $B=\{s=0\}$ the local model is $u^2=s$. At each node the monodromy along one of the loops
fixes the first four digits and changes $e$. For example  $(y_p\delta_p)^2$ over $(p,p)$ for $p=1,2$, and $x_4y_3$ over $(4,3)$.
It follows that $S$ is smooth, and $f$ is branched exactly along $B$ and at the six nodes.

\begin{proposition}\label{smooth}
$S$ is a fake quadric with $H_1(S;\ZZ)\cong(\ZZ/2)^3\oplus\ZZ/4$.
\end{proposition}

\begin{proof}
  Clearly $e(S)=4$. The computation of $H_1(S;\ZZ)$ is straightfoward from the monodromy description. As a corollary we have
  $b_1(S)=0$ and $b_2(S)=2$. The preimages
  $F_1$ and $F_2$ of a general
horizontal and a general vertical line are double covers of $C_2$ and
$C_1$ branched at $0$ and $8$ points. It follows that $b_2^+=1$ and we also see that  $S$ is minimal surface of 
of general type.
\end{proof}

\section{A projective plane in $S$}\label{sec:sphere}

Let $\gamma$ be the upper half of the circle of radius $1$ about $3$,
including its endpoints $2$ and $4$, and let $\Delta\subset\gamma\times\gamma$ be the triangle bounded by $z=2$,
$w=4$ and the diagonal $z=w$.

Its interior
misses $D$. Move the base point into $\Delta$, first $z$ from $i$ to
$3+i$, then $w$ from $-i$ to $\frac72-i$ and straight up to $\gamma$. The meridians of $z=2$, of the diagonal and of $w=4$ are then $x_2$,
$\delta_4$ and $y_4$. On the sheets with $a=b=0$, written $(c,d,e)$,
$$
000\xrightarrow{x_2}100\xrightarrow{\delta_4}110
\xrightarrow{x_2}011\xrightarrow{\delta_4}001
\xrightarrow{x_2}101\xrightarrow{\delta_4}111
\xrightarrow{x_2}010\xrightarrow{\delta_4}000.
$$
The copies of $\Delta$ on these eight sheets form an octagon around a point over $(2,2)$, consecutive
copies sharing their side over $z=2$ or over the diagonal.
Opposite sheets differ by $\tau=y_4^2$, so opposite copies share their side over $w=4$,
and identifying opposite
sides gives a map $\varphi\colon\RP^2\to S$.

Let $S'\to S$ be the double cover obtained by appending a digit $h$,
where $x_3,x_4$ add $c$ to $h$, $y_1,y_4$ add $b$, and the other generators fix $h$. The relations still hold,
and every meridian has the same cycle lengths on $S'$ as on $S$,
so $S'\to S$ is unramified. Since $x_2$ and $\delta_4$ fix $h$ while $y_4^2$ changes it, the
octagon lifts to $S'$. However  its opposite sides no longer match. So the
loop in $\RP^2$ through the octagon from a side to the opposite side
does not lift to a closed loop in $S'$, and $\varphi$ is nontrivial on $\pi_1$.

So $\pi_1(S)$ has an element of order two. Since the fundamental group of a closed aspherical manifold is torsion-free, $S$ is not
aspherical, and its universal cover is not $\HH\times\HH$.

\begin{remark}
One can check from the local models at the vertices that $\varphi$ is an embedding onto a smooth, totally real projective plane. Its
preimage in $S'$ is then a totally real sphere with
self-intersection $-2$, which gives an explicit nonzero element of $\pi_2(S')\cong\pi_2(S)$.
\end{remark}

\end{document}